\documentclass[12pt]{article}
\usepackage{amsmath,amssymb}
\usepackage{fullpage}

\newtheorem{thm}{Theorem}[section]
\newtheorem{prop}[thm]{Proposition}
\newtheorem{lemma}[thm]{Lemma}
\newtheorem{cor}[thm]{Corollary}

\newcommand{\N}{{\mathbb N}}
\newcommand{\E}{{\mathbb E}}

\newcommand{\cL}{{\cal L}}
\newcommand{\cM}{{\cal M}}
\newcommand{\approxsb}[1]{\mathrel{\mathop{\approx}\limits_{#1}}}
\newcommand{\modo}[1]{{\left|#1\right|}}
\newcommand{\normo}[1]{{\left\|#1\right\|}}
\newcommand{\smodo}[1]{{\mathopen|#1\mathclose|}}
\newcommand{\snormo}[1]{{\mathopen\|#1\mathclose\|}}
\newcommand{\Proof}{\medskip\noindent{\bf Proof:}\ \ }
\newcommand{\Proofof}[1]{\medskip\noindent{\bf Proof of #1:}\ \ }
\newcommand{\qed}{\leavevmode\unskip\penalty9999
                  \hbox{}\nobreak\hfill$\square$\goodbreak\medskip}

\begin{document}
\title{Measuring the magnitude of sums\\
of independent random variables}

\author{%
Pawe{\l} Hitczenko%
\thanks{The first named author was partially supported 
by NSF grant DMS 9401345.}\\
Department of Mathematics\\
North Carolina State University\\
Raleigh, NC 27695-8205\\
pawel@math.ncsu.edu\\
http://www.math.ncsu.edu/\~{}pawel\\
\\
Stephen Montgomery-Smith%
\thanks{The second named author was 
partially supported
by NSF grants DMS 9424396 and DMS 9870026, and by the
University of Missouri Research Board.}\\
Department of Mathematics\\
University of Missouri-Columbia\\
Columbia, MO 65211\\
stephen@math.missouri.edu\\
http://math.missouri.edu/\~{}stephen}

\date{}

\maketitle

\begin{abstract}
\noindent
This paper considers how to measure the magnitude of the sum of 
independent random variables in several ways.
We give a formula for the tail distribution for
sequences that satisfy the so called L\`evy property.
We then give a connection between the tail distribution and the
$p$th moment, and between the $p$th moment and
the rearrangement invariant norms.
\\
\\
Keywords: sum of independent random variables, 
tail distributions, decreasing rearrangement, 
$p$th moment, rearrangement invariant space,
disjoint sum, maximal function,
Hoffmann-J{\o}rgensen/Klass-Nowicki Inequality,
L\`evy Property.
\\
\\
A.M.S.\ Classification (1991): Primary 60G50, 60E15, 46E30; 
Secondary 46B09. 

\end{abstract}

\section{Introduction}
\label{intro}

This paper is about the following type of problem: 
given independent (not necessarily identically distributed) random
variables 
$X_1$, $X_2,\dots,$ $ X_N$, find the `size' of $\smodo S$, 
where 
$$ S =\sum_{n=1}^N X_n .$$ 

We will examine several ways to measure this size.  The first will be through
tail distributions, that is, $\Pr(\smodo S > t)$.
Finding an exact
solution to this problem would be a dream of probabilists,  so we
have to 
temper our desires in some manner. In fact, this problem goes back to
the 
foundations of probability in the following form: if the sequence
$(X_n)$ 
consists of random variables that are mean zero, identically
distributed and have finite variance, find the asymptotic value of 
$\Pr(\smodo S > \sqrt N t)$ as $N \to \infty$. This is answered, of
course, 
by the Central Limit Theorem, which tells us that the answer is the
Gaussian 
distribution. There has been a tremendous amount of work on
generalizing 
this. We refer the reader to almost any advanced work on probability. 

Our approach is different. Instead of seeking asymptotic solutions,
we will look for approximate solutions. That is, we seek a function
$f(t)$, 
computed  from $(X_n)$, such that there is a positive constant $c$ with 
$$ c^{-1} f(c t) \le \Pr(|S|>t) \le c \, f(c^{-1} t) .$$ 

The second 
measurement of the size of $\smodo S$ will be through
the $p$th moments, $\snormo S_p$.  
Again, we shall be searching for approximate solutions, that is, finding
a quantity $A$ such that there is a positive constant $c$ so that
$$ c^{-1} A \le \snormo S_p \le c \, A .$$
While this may seem like quite a
different problem, in fact, as we will show, there is a precise connection
between the two, in that obtaining an approximate formula for $\snormo S_p$
with constants that are uniform as $p\to \infty$ is equivalent to obtaining
an approximate formula for the tail distribution.

The third way that we shall look at is to find the size of $\smodo S$ in
a rearrangement invariant space.  This line of research was began by
Carothers and Dilworth (1988) who obtained results for Lorentz spaces,
and was completed by Johnson and Schechtman (1989).  Our results will
give a comparison of the size of $\smodo S$ in the rearrangement invariant
space with $\snormo S_p$, obtaining a greater control on the sizes of the
constants involved than the previous works.

Many of the results of this paper will be true of all sums of
independent random variables, even those that are vector valued, with
the following proviso.  Instead of considering 
the sum $S = \sum_n X_n$, we will
consider the maximal function $U = \sup_n \modo{\sum_{k=1}^n X_k}$.  
We will define a property for sequences called the {\em L\`evy property},
which will imply that $U$ is comparable to $S$.  Sequences
with this L\`evy property will include positive random
variables, symmetric random variables, and identically distributed
random variables. The result of this paper that gives the tail
distribution for $S$ is only valid for real valued sequences of
random variables that satisfy the L\`evy property.  However the results
connecting the $L_p$ and the rearrangement invariant norms to the
tail distributions of $U$ are valid for
all sequences of vector valued independent random variables.

Let us first give the historical context for these results, considering
first the problem of approximate formulae for the tail distribution.
Perhaps the earliest works are the
Paley-Zygmund inequality (see for example Kahane (1968, Theorem 3, Chapter 2)), 
and Kolmogorov's reverse maximal inequality (see for example Shiryaev (1980,
Chapter 4, 
section 2.)) Both give (under an extra assumption) a lower bound on the
probability that a sum of independent, mean zero random variables
exceeds a 
fraction of its standard deviation and both may be regarded as a sort
of 
converse to the Chebyshev's inequality.
Next, in 1929, Kolmogorov, proved a
two-sided exponential inequality for sums of independent, mean-zero, 
uniformly bounded, random variables (see for example Stout
(1974, 
Theorem 5.2.2) or Ledoux and Talagrand (1991, Lemma 8.1)).
All of these results require some restriction on the nature
of the sequence $(X_n)$, and on the size of the level $t$.

Hahn and Klass (1997)
obtained very good bounds on one sided tail probabilities for sums
of independent, identically distributed, real valued 
random variables.  Their
result had no restrictions on the nature of the random variable, or
on the size of the level $t$.
In effect, their
result worked by removing the very large parts of the random variables,
and then using an exponential estimate on the rest.  
We will take a similar approach in this paper.

Let us next look at the $p$th moments.  Khintchine (1923) gave
an inequality for Rademacher (Bernoulli) sums.  
This very important 
formula has found extensive applications in analysis and probability.
Khintchine's result was extended
to any sequence of positive or mean zero random variables by 
the celebrated result of Rosenthal (1970).
The order of the best constants as $p \to \infty$ was 
obtained by Johnson, Schechtman and Zinn (1983),
and Pinelis (1994) refined this still further.  Now even more precise
results are known, and we refer the reader to 
Figiel, Hitczenko, Johnson, Schechtman and Zinn (1997).  However, the
problem with all these results is that the constants were not uniformly
bounded as $p \to \infty$.

Khintchine's inequality was generalized independently by 
Montgomery and Odlyzko (1988) and Montgomery-Smith (1990).  They
were able to give approximate bounds on the tail probability for 
Rademacher sums, with no restriction on the level $t$.
Hitczenko (1993) obtained an approximate formula for the
$L_p$ norm of Rademacher sums, where the constants were uniformly bounded
as $p \to \infty$.  (A more 
precise version of this last result was obtained in Hitczenko-Kwapie\'n
(1994) and it was used to give a simple proof of the lower bound
in Kolmogorov's exponential inequality.) 

Continuing in the direction of Montgomery and Odlyzko, Montgomery-Smith and
Hitczenko,
Gluskin and Kwapie\'n (1995) extended tail and moment estimates
from Rademacher sums to weighted sums of random variables with
logarithmically concave tails (that is, $P(\modo{X}\ge t)=\exp(-\phi(t))$,
where 
$\phi:[0,\infty)\to[0,\infty)$ is convex).  
After that, Hitczenko, Montgomery-Smith, 
and Oleszkiewicz (1997) treated the case of logarithmically
convex tails 
(that is, the $\phi$ above is concave rather than convex). 
It should be 
emphasized that in the last paper, the result of
Hahn and Klass (1997) played a critical role.

The breakrough came with the paper of Lata{\l}a (1997), who solved
the problem of finding upper and lower bounds for general
sums of positive or symmetric random variables, with uniform constants
as $p \to \infty$.
His method made beautiful use of special properties of the function 
$t \mapsto t^p$. 
In a short note, Hitczenko and 
Montgomery-Smith (1999) showed how to use Lata{\l}a's result to derive 
upper and lower bounds on 
tail probabilities.  Lata{\l}a's result is the primary
motivation for this paper.

The main tool we will use is the Hoffmann-J{\o}rgensen Inequality. In
fact, we will use a stronger form of this inequality, due to
Klass and Nowicki (1998).
The principle in many of our proofs is the following idea.  Given 
a sequence of random variables $(X_n)$, we choose an appropriate level 
$s>0$.  Each random variable $X_n$ is split into the 
sum $X^{(\le s)}_n + X^{(>s)}_n$, where $X^{(\le s)}_n = X_n
I_{\modo{X_n}\le 
s}$, and $X^{(>s)}_n = X_n I_{\modo{X_n} > s}$.  It turns out that the
quantity $(X_n^{(>s)})$ can either be disregarded, or it can be  
considered as a sequence of disjoint random variables. 
(By ``disjoint'' we mean that the random variables are disjointly
supported as functions on the underlying probability space.) As for the
quantity 
$\sum_n X_n^{(\le s)}$, it will turn out that the level $s$ allows one
to apply the Hoffmann-J{\o}rgensen/Klass-Nowicki Inequality so that it may 
be compared with
quantities that we understand rather better.

Let us give an outline of this paper.  In Section~\ref{Defins}, we will give
definitions. This will include the notion of decreasing rearrangement,
that is, the inverse to the distribution function.  Many results of
this paper will be written in terms of the decreasing rearrangement.
Section~\ref{Klass-nowicki} is devoted to the Klass-Nowicki Inequality.  
Since our result is
slightly stronger than that currently in the literature, we will
include a full proof. 
In Section~\ref{Levy}, we will introduce and discuss the
L\`evy property.  This will include a ``reduced comparison principle''
for sequences with this property. 
Section~\ref{Tail} contains the formula for the tail distribution of sums of
real
valued random variables. 
Then in Section~\ref{Lp}, we demonstrate the connection between $L_p$-norms 
of such sums and their tail distributions.
In Section~\ref{Ri} we will discuss sums of independent random variables in
rearrangement invariant spaces.

\section{Notation and definitions}
\label{Defins}

Throughout this paper, a random variable will be a measurable function
from a probability space to some Banach space (often the real
line). The norm in the implicit Banach space will always be denoted by
$\smodo{\ \cdot\ }$.

Suppose that $f:[0,\infty) \to [0,\infty]$ is a decreasing
function.  Define the
{\em left continuous inverse\/} to be
$$ f^{-1}(x-) = \sup\{ y: f(y) \ge x \} ,$$
and the {\em right continuous inverse\/} to be
$$ f^{-1}(x+) = \sup\{ y: f(y) > x \} .$$

In describing the tail distribution of a random variable $X$, instead
of considering 
the function $t \mapsto \Pr(\smodo X > t)$, we will 
consider its right continuous inverse, which we will denote by $X^*(t)$. 
In fact, this
quantity appears very much in the literature, and is more commonly
referred to as the {\em decreasing rearrangement\/} of $\smodo
X$. Notice that if one considers $X^*$ to be a random variable on the
probability space $[0,1]$, then $X^*$ has exactly the same law as
$\smodo X$.
We might also consider the left continuous inverse $t \mapsto
X^*(t-)$. Notice that $X^*(t) \le x \le X^*(t-)$ if and only if $\Pr(\smodo
X > x) \le t \le \Pr(\smodo X \ge x)$.

If $A$ and $B$ are two quantities (that may depend upon certain 
parameters),  we will write $A \approx B$ to mean that there exists
positive constants $c_1$ and $c_2$ such that $c_1^{-1} A \le B \le c_2
A$. We will call $c_1$ and $c_2$ the constants of approximation. If
$f(t)$ and $g(t)$ are two (usually decreasing) functions on
$[0,\infty)$, we will write $f(t) \approxsb t g(t)$ if there exist
positive constants $c_1$, $c_2$, $c_3$ and $c_4$ such that $c_1^{-1}
f(c_2 t) \le g(t) \le c_3 f(c_4^{-1} t)$\ for all $t \ge 0$. Again, we
will call $c_1$, $c_2$, $c_3$ and $c_4$ the constants of
approximation.

Suppose that $X$ and $Y$ are random variables. Then the statement 
$\Pr(\smodo X > t) \approxsb t \Pr(\smodo Y > t)$ is the same as the
statement $X^*(t) \approxsb t Y^*(t)$.  Since $X^*(t) = 0$ for $t \ge 1$
the latter statement is equivalent to
the existence of positive constants $c_1$, $c_2$, $c_3$, $c_4$ and
$c_5$ such that $c_1^{-1} X^*(c_2 t) \le Y^*(t) \le c_3 X^*(c_4^{-1}
t)$ 
for $t \le c_5^{-1}$. 

To avoid bothersome convergence problems, we will always suppose that
our sequence of independent random variables $(X_n)$ is of finite
length. Given a sequence of independent random variables $(X_n)$, when
no 
confusion will arise, we will use the following notations. If $A$ is a
finite 
subset of $\N$, we will let $ S_A = \sum_{n \in A} X_n$,  and $M_A = 
\sup_{n \in A} \smodo{X_n}$. If $k$ is a positive integer, then 
$S_k = S_{\{1,\dots,k\}}$ and $M_k = M_{\{1,\dots,k\}}$. We will 
define the maximal function
$U_k = \sup_{1 \le n \le k}\modo{S_n}$. Furthermore, $S = S_N$, $M =
M_N$, and $U = U_N$, where $N$ is the length of the sequence $(X_n)$.

If $s$ is a real number, we will write 
$X_n^{(>s)} = X_n I_{\smodo{X_n} > s} $ 
and $X_n^{(\le s)} = X_n I_{\smodo{X_n} \le s} = X_n - X_n^{(>s)}$.
For $A \subset\N$, we will write $S^{(\le s)}_A = \sum_{n
\in A} X^{(\le s)}_n$. 
Similarly we define $S^{(>s)}_A$, 
 $S^{(\le s)}_k$, etc.

Another quantity that we shall care about is the decreasing
rearrangement of the disjoint sum of random variables.  
This notion was used by Johnson, Maurey,
Schechtman  and Tzafriri (1979), Carothers and Dilworth (1988),
and Johnson and Schechtman (1989), all in the context of sums of
independent random variables.
The disjoint
sum of the sequence $(X_n)$ is the measurable function on the measure
space $\Omega \times \N$ that takes $(\omega,n)$ to
$X_n(\omega)$.  We shall denote the decreasing rearrangement of
the disjoint sum by $\tilde \ell:[0,\infty) \to [0,\infty]$, that is,
$\tilde \ell(t)$ is the least number such that 
$$ \sum_n \Pr(\smodo{X_n} > \tilde \ell(t)) \le t .$$
Define $\ell(t)$ to be $\tilde \ell(t)$ if $0 \le t \le
1$, and $0$ otherwise. Since $\ell(t)$ is only non-zero when $0 \le t
\le 1$, we will think of $\ell$
as being a random variable on the 
probability space $[0,1]$ with Lebesgue measure.
The quantity $\ell$ is effectively $M$ in disguise.  This next result 
essentially appears in Johnson and Schechtman (1989).

\begin{prop}
\label{ell-max}
 If $0<t < 1$, then 
$$ \ell(2t) \le \ell(t/(1-t)) \le M^*(t) \le\ell(t) .$$
\end{prop}

\Proof  The first inequality follows easily
once one notices that both sides of this inequality are zero if $t >
1/2$. 

To get the second inequality, note that, by an easy argument, if $\alpha_1$, 
$\alpha_2,\dots \ge 0$ with $\sum_n \alpha_n \le 1$, then 
$$ 1-\sum_n \alpha_n \le \prod_n (1-\alpha_n)\le 1 - {\frac{\sum_n \alpha_n
}{ 1 
+ \sum_n \alpha_n}} .$$ So, if $\Pr(\ell > x) = \sum_n \Pr(\smodo{X_n}
> x) 
\le 1$, then $$ \Pr(M > x) = 1-\prod_n (1-\Pr(\smodo{X_n} > x)) ,$$  
and hence $$ {\frac{\Pr(\ell > x) }{ 1+\Pr(\ell>x)}} \le \Pr(M>x) \le 
\Pr(\ell > x) .$$ Taking inverses, the result follows. 
\qed

\section{The Klass-Nowicki Inequality}
\label{Klass-nowicki}

This section is devoted to the following result --- the Klass-Nowicki
Inequality.

\begin{thm}
\label{klass-nowicki}
Let $(X_n)$ be a sequence of independent random
variables. Then for all positive integers $K$ we have 
$$ \Pr(U > (3 K - 1) t)\le {\frac{1 }{ K!}} \left(\frac{\Pr(U > t) }{ 1 -
\Pr(U > t) }\right)^{K}+ \Pr(M > t) .$$
\end{thm}

The original inequality of this form was 
for Rademacher (or Bernoulli) sums and $K=2$, and was due to Kahane
(1968). This was extended by Hoffmann-J{\o}rgensen (1974)
to general sums, at
least for positive or symmetric random variables, for the case $K =
2$.  
Indeed, if one wants Theorem~\ref{klass-nowicki}
for $K>2$, but without 
the $K!$ factor, this may be obtained by iterating the
Hoffmann-J{\o}rgensen Inequality, as was done by Johnson and 
Schechtman (1989, Lemmas 6 and 7).
(Both Kahane and Hoffmann-J{\o}rgensen obtained slightly different
constants 
than those we have presented.  Also, in neither case 
did a factor like $(1-\Pr(U>t))$ 
appear in their formulae.) 

Klass and Nowicki (1998) were able to obtain Theorem~\ref{klass-nowicki}, 
at least in the case when the random variables are positive or
symmetric. (However their constants are better than ours.) Removing
the 
positive or symmetric condition is really not so hard, but because 
it does not appear in the literature in this manner, we will give a 
complete proof of Theorem~\ref{klass-nowicki}.  

We also note that this inequality has some comparison with a result
that 
appears in Ledoux and Talagrand (1991, Theorem 6.17.)

\Proof  Let $N$ be the length of the sequence $(X_n)$. During this
proof, let us write $(m,n]$ for the set of integers greater than $m$
and not greater than $n$.

We start with the observation 
$$
  \Pr(U >  (3 K - 1) t) 
  \le 
  \Pr(U >(3 K - 1) t \text{ and } M \le t) + \Pr(M > t) .
$$
Now, if we have that both $U > (3K-1)t$ and $M \le t$, then we ensure
the existence of an increasing sequence of non-negative integers
$m_0,\dots,$\ 
$m_{K}$, bounded by $N$, and defined as follows. Set $m_0 = 0$. 
If we have picked $m_{l-1}$, let $m_l$ be the smallest positive
integer greater than $m_{l-1}$ such that $\smodo{S_{(m_{l-1},m_l]}} >
2t$.  For $l = 1$, it is clear that such an integer exists.
Let us explain why the  integer $m_l \le N$ exists if $2 \le l \le K$.

For $1 \le l' \le l-1$, and $m_{l'-1} < k \le m_{l'}-1$, we have that
$\smodo{S_{(m_{l'-1},k]}} \le 2t$ and $\smodo{X_{m_{l'}}} \le t$.  Hence
$$ \smodo{S_{m_{l'}}} 
   =
   \modo{ \sum_{j=1}^{l'} S_{(m_{j-1},m_j-1]} + X_{m_j} }
   \le
   3 l' t ,$$
and
$$ \smodo{S_{k}} 
   =
   \modo{ \left(\sum_{j=1}^{l'-1} S_{(m_{j-1},m_j-1]} + X_{m_j}\right)
          + S_{(m_{l'-1},k]} }
   \le
   (3 l' - 1) t .$$
But we know that there exists a number $m$ such that
$\smodo{S_m} > (3K-1)t$. Hence, we must 
have that $m > m_{l-1}$, and that 
$\smodo{S_{(m_{l-1},m]}} > (3K-1)t - 3(l-1)t \ge 2t$.

Therefore 
$$ \Pr( U > (3K-1) t \text{ and } M \le t) \le \sum_{1 \le m_1
< \cdots < m_{K} \le N}p_{0,m_1} p_{m_1,m_2} \cdots p_{m_{K-1},m_K}
,$$ where $$ p_{m,n} = \Pr(\text{$\smodo{S_{(m,k]}} \le 2t$ for $m \le
k < n$, and $\smodo{S_{(m,n]}} > 2t$}) .$$ 
Now let us show the following inequality: 
$$ \sum_{k=m+1}^n p_{m,k}
   \le
   {\frac{1}{ 1 - \Pr(U>t)}}\sum_{k=m+1}^n \tilde p_k ,$$ 
where 
$$ \tilde p_n = \Pr(\text{$\smodo{S_k} \le t$ for $1 \le k < n$, and 
   $\smodo{S_n} > t$}) .$$ 
Using independence, we have that 
\begin{eqnarray*}
  \sum_{k=m+1}^n p_{m,k}
  &=&\Pr(\sup_{m < k \le n} \smodo{S_k - S_m} > 2t) \\ 
  &=&\Pr(\sup_{m < k \le n} \smodo{S_k - S_m} > 2t \big| U_m \le t) \\ 
  &\le& \Pr(\sup_{m < k \le n} \smodo{S_k} > t \big| U_m \le t) \\
  &=& {\frac{\Pr(\sup_{m<k\le n} \smodo{S_k} > t 
      \text{ and }\sup_{1 \le k \le m} \smodo{S_k} \le t)
      }{ \Pr(U_m \le t) }} \\ 
  &\le& {\frac{1}{ 1 - \! \Pr(U>t)}}\sum_{k=m+1}^n \tilde p_k ,
\end{eqnarray*} 
as required.

Now we rearrange the sum as follows: 
\begin{eqnarray*}
  && \sum_{1 \le m_1 < \cdots < m_{K} \le N}p_{0,m_1} p_{m_1,m_2} 
    \cdots p_{m_{K-1},m_K}\\
  &=& \sum_{1 \le m_1 < \cdots < m_{K-1} \le N}p_{0,m_1} p_{m_1,m_2} 
      \cdots p_{m_{K-2},m_{K-1}} \sum_{m_K = m_{K-1}+1}^Np_{m_{K-1},m_K} \\
  &\le& {\frac{1 }{ 1 - \Pr(U>t)}}\sum_{1 \le m_1 < \cdots < m_{K-1} \le N}
        p_{0,m_1} p_{m_1,m_2} \cdots p_{m_{K-2},m_{K-1}} 
       \sum_{m_K = m_{K-1}+1}^N \tilde p_{m_K} .
\end{eqnarray*} 
Now we rearrange this last quantity to get 
\begin{eqnarray*}
  &&\sum_{1 \le m_1 < \cdots < m_{K} \le N} p_{0,m_1} p_{m_1,m_2} 
    \cdots p_{m_{K-1},m_K}\\
  &\le& {\frac{1 }{ 1 - \Pr(U>t)}}
        \sum_{1 \le m_1 < \cdots < m_{K-2} < m_K\le N} 
        p_{0,m_1} p_{m_1,m_2} \cdots p_{m_{K-3},m_{K-2}} \tilde p_{m_K} \times\\
  && \times
        \sum_{m_{K-1} = m_{K-2}+1}^{m_K} p_{m_{K-2},m_{K-1}} \\ 
  &\le& {\frac{1 }{ (1- \Pr(U>t))^2}}
        \sum_{1 \le m_1 < \cdots < m_{K-2} < m_K\le N}
        p_{0,m_1} p_{m_1,m_2}
        \cdots p_{m_{K-3},m_{K-2}} \tilde p_{m_K} \times\\
  && \times
        \sum_{m_{K-1} = m_{K-2}+1}^{m_K} \tilde p_{m_{K-1}} .
\end{eqnarray*}
Repeating this 
argument $(K-2)$ more times, we eventually see that 
\begin{eqnarray*}
  && \sum_{1 \le m_ < \cdots < m_K \le N} p_{0,m_1} p_{m_1,m_2} \cdots 
     p_{m_{K-1},m_K}\cr 
  &\le& {\frac{1 }{ (1 - \Pr(U>t))^K}}\sum_{1 \le m_1 < \cdots < m_K \le N}
        \tilde p_{m_1} \tilde p_{m_2} \cdots \tilde p_{m_K} .
\end{eqnarray*}
Now, since $K$ distinct numbers may be rearranged in $K!$ different ways,
we have that 
\begin{eqnarray*}
  && \sum_{1 \le m_1 < \cdots < m_K \le N} 
     \tilde p_{m_1} \tilde p_{m_2} \cdots \tilde p_{m_K}\\ 
  &=& {\frac{1 }{ K!}}
      \sum_{\substack{1 \le m_1,m_2,\dots,m_K \le N\\ m_1,m_2,\dots,m_K 
      \text{ distinct}}} \tilde p_{m_1} \tilde p_{m_2} \cdots \tilde
      p_{m_K}\\ 
  &\le& {\frac{1 }{ K!}}
        \sum_{1 \le m_1,m_2,\dots,m_K \le N} \tilde p_{m_1} \tilde p_{m_2}
        \cdots \tilde p_{m_K}\\
  &=& {\frac{1 }{ K!}} \left( \sum_{k=1}^N \tilde p_{k} \right)^K . 
\end{eqnarray*} 
Since $$ \sum_{k=1}^N \tilde p_k= \Pr(U > t),$$ we obtain the
result.
\qed

Let us now understand what this result means in terms of the
decreasing rearrangement.

\begin{cor}
\label{klass-nowicki-decr}
There exists  a universal positive constant 
$c_1$ such that for any sequence of independent random variables
$(X_n)$, and for $0<t\le s \le 1/2$ we 
have $$ U^*(t) \le c_1 {\frac{\log(1/t) }{
\max\{\log(1/s),\log\log(4/t)\}}} \bigl(U^*(s) + M^*(t/2)\bigr). $$
\end{cor}

\Proof  Notice that if $f,g:[0,\infty)\to[0,\infty]$ are decreasing functions,
 then $(\max\{f,g\})^{-1} = \max\{f^{-1},g^{-1}\}$, and if $f \le g$, 
then $f^{-1} \le g^{-1}$, where here $f^{-1}$ denotes either the left
or right continuous inverse of $f$. Since $A+B \le 2 \max\{A,B\}$ for
any 
two positive numbers $A$ and $B$, from Theorem~\ref{klass-nowicki}, we have 
that if 
$\Pr(U>t) \le 1/2$, then for all positive integers $K$ 
$$ \Pr(U > (3 K - 1) t)\le 2 \max \left\{{\frac{1 }{ K!}} 
\bigl(2\Pr(U > t)\bigr)^K, \Pr(M > t) \right\} .$$ Taking inverses, we
see 
that if $(K!t/2)^{1/K} \le 1/2$, then $$ {\frac{1 }{3K-1}} U^*(t) \le
\max\left
\{ U^*\left({\frac12}\left({\frac{K!t }{ 2}}\right)^{1/K}\right) ,
M^*\left({\frac{t}{2}}
\right) \right\} .$$ 
Now, using the fact that $\max\{A,B\} \le A+B$
for any 
positive numbers $A$ and $B$, and by choosing $K$ to be the smallest
integer such that $s \le(K!t / 2)^{1/K}$, and by some elementary but
tedious algebra, the result follows. \qed

\section{The L\`evy Property}
\label{Levy}

Let $(X_n)$ be a sequence of independent random variables. We will say 
that $(X_n)$ satisfies the {\em L\`evy property\/} with constants $c_1$ 
and $c_2$ if whenever $A
\subseteq B 
\subseteq \N$, with $A$ and $B$ finite, then $$
\Pr(\smodo{S_A} > c_1 t)\le c_2 \Pr(\smodo{S_B} > t) .$$ 
The casual reader should beware that this property has nothing to do
with L\`evy processes.

The sequence $(X_n)$ has the {\em strong L\`evy property\/} with
constants 
$c_1$ and $c_2$ if for all $s>0$ the sequence $(X_n^{(\le s)})$ has
the L\`evy property with constants $c_1$ and $c_2$.

Here are examples of sequences with the strong L\`evy property. (It
may be easily seen that in all these cases that it is sufficient to
show that 
they have the L\`evy property.) 
\begin{enumerate}
\item All positive sequences,
with constants $1$ and $1$. 
\item All symmetric sequences (that
is, when $X_n$ has the same law as $-X_n$) with constants $1$ and
$2$. This ``reflection property'' plays a major role in results
attributed to L\`evy, hence the name of the property. 
\item Sequences of identically distributed random variables. This was shown 
independently by Montgomery-Smith (1993) with constants $10$ and $3$,
and by Lata{\l}a (1993) with constants $5$ and $4$, or $7$ and
$2$.
\end{enumerate}
We see that sequences with the L\`evy property
satisfy a maximal inequality.

\begin{prop}
\label{max}
Let
$(X_n)$ be a sequence of independent random satisfying the L\`evy
property with constants $c_1$ and $c_2$. Then $$ \Pr( U > 3 c_1 t
) \le3 c_2 \Pr ( \smodo S > t ) .$$ Thus $M^*(t) \le 6 c_1
S^*(t/3c_2)$.
\end{prop}

\Proof
The first statement is an immediate corollary of the
following result known as L\`evy-Ottaviani
inequality: 
$$ \Pr(U_N > 3 t)\le3 \sup_{1 \le k \le N} \Pr( \smodo{S_k} > t ) .$$
(Billingsley (1995, Theorem 22.5, p. 288) attributes
this result to Etemadi (1985) who proved it with constants 4 in both
places, but the same proof gives constants 3; see, for example,
Billingsley.  
However the first named author learned this result from
Kwapie\'n in 1980.)

The second statement follows from the first, since $M
\le 2 U$. 
\qed

We end with a lemma that lists some elementary properties.  
Part~\ref{levy-comparison}
of the lemma might be thought of as a kind of reduced 
comparison principle.

\begin{lemma}
\label{levy}
Let $(X_n)$ be a sequence of random variables 
satisfying the strong L\`evy property. 
\begin{enumerate}
\item 
\label{levy-comparison}
There exist  positive constants $c_1$ and $c_2$, depending
only upon the L\`evy constants of $(X_n)$, such that if $s \le 1/2$
and 
$0 \le t \le 1$, then $$ (S^{(\le M^*(s))})^*(t) \le c_1
S^*(c_2^{-1} t) .$$ 
\item 
\label{levy-incsuper}
There exist positive constants $c_1$
and $c_2$, depending only upon the 
strong L\`evy constants of $(X_n)$, such that if $r \le s \le 1/2$,
and if $0 \le t \le 1$, then
$(S^{(\le M^*(s))})^*(t) \le c_1 (S^{(\le M^*(r))})^*(c_2^{-1}t)$.
\item 
\label{levy-decsuper}
If $0 \le s \le t \le 1$, then $S^*(t) \le (S^{(\le 
M^*(s))})^*(t-s)$, and
 $(S^{(\le M^*(s))})^*(t) \le S^
*(t-s)$. In particular, $S^*(t) \le (S^{(\le 
M^*(t/2))})^*(t/2)$, and $(S^{(\le 
M^*(t/2))})^*(t) \le S^*(t/2)$. 
\item 
\label{levy-equiv}
For $\alpha,\beta > 0$, we have
that 
$$ (S^{(\le M^*(t))})^*(t)
\approxsb t (S^{(\le M^*(\alpha t))})^*(\beta t) $$
where the 
constants of approximation depend only upon $\alpha$, $\beta$ and the strong
L\`evy constants of $(X_n)$. 
\item 
\label{levy-ellequiv}
We have that $$ S^*(t) \approxsb t (S^{(\le M^*(t))})^*(t) \approxsb t
(S^{(\le \ell(t))})^*(t) ,$$ where 
the constants of approximation depend only upon the strong L\`evy
constants of $(X_n)$.
\end{enumerate}
\end{lemma}

\Proof  
Let us start with part~\ref{levy-comparison}.
For each set $A \subseteq \N$, define the event
$$ E_A = \{ \smodo{X_n} \le M^*(s) \text{ if and only if } n \in A \} .$$
Note that the whole probability space is the disjoint union of these
events.  Also
$$ \{ \smodo{S^{(\le M^*(s))}} > x \} \cap E_A
   =
   \{ \smodo{S_A} > x \} \cap E_A .$$
Furthermore, by independence, we see that
\begin{eqnarray*}
   &&\Pr(\smodo{S_A} > x \text{ and } E_A) \\
   &&= \Pr(\smodo{S_A} > x \text{ and } \smodo{X_n} \le M^*(s) 
       \text{ for } n \in A)
       \Pr(\smodo{X_n} > M^*(s) \text{ for } n \notin A) .
\end{eqnarray*}
Hence
\begin{eqnarray*}
  &&\Pr(\smodo{S^{(\le M^*(s))}} > x) \\
  &&= \sum_{A \subseteq \N}
      \Pr(\smodo{S_A} > x \text{ and } \smodo{X_n} \le M^*(s) 
      \text{ for } n \in A)
      \Pr(\smodo{X_n} > M^*(s) \text{ for } n \notin A) \\ 
  &&\le
      2 \sum_{A \subseteq \N}
     \Pr(\smodo{S_A} > x)
     \Pr(\smodo{X_n} \le M^*(s) \text{ for } n \in A) \Pr(\smodo{X_n} >
     M^*(s)
     \text{ for } n \notin A) \\
  &&\le c_2 \Pr(\smodo S > c_1^{-1} x), 
\end{eqnarray*}
where here we have used the fact that 
$$ \Pr(\smodo{X_n} \le M^*(s) \text{ for } n \in A) \ge \Pr(M \le M^*(s)) 
\ge 1-s \ge 1/2 .$$ 
Part~\ref{levy-incsuper} follows by applying 
part~\ref{levy-comparison} to $S^{(\le M^*(r))}$.

Part~\ref{levy-decsuper} follows from the observation that $$ 
\Pr(S \ne S^{(\le M^*(s))})\le\Pr(M > M^*(s)) \le s .$$ 
Hence, if $\Pr(S > \alpha) \ge t$, then $\Pr(S^{(\le M^*(s))} >
\alpha) 
\ge t-s$, and conversely, if $\Pr(S^{(\le M^*(s))} > \alpha) > t$ then
$\Pr(S >\alpha) \ge t-s$.

To show 
part~\ref{levy-equiv},
we may suppose without loss of generality that $\alpha=1$ and $\beta > 1$.
Clearly $S^{(\le M^*(t))}(t) \ge S^{(\le M^*(t))}(\beta t)$, so we need
only show an opposite inequality.  From part~\ref{levy-incsuper}, there are
positive
constants $c_1$ and $c_2$, depending only upon the strong L\`evy constants of
$(X_n)$, such that for $0\le t \le 1/2$
$$  S^{(\le M^*(t))}(t) 
    \le c_1 S^{(\le M^*(c_2^{-1}\beta^{-1})t)}(c_2^{-1} t) 
    \le c_1 S^{(\le M^*(c_3^{-1}t)}(c_3^{-1} \beta t) ,$$
where $c_3 = c_2 \beta$.

Part~\ref{levy-ellequiv} follows easily by combining
part~\ref{levy-decsuper}, part~\ref{levy-equiv},
and Proposition~\ref{ell-max}. \qed

\section{Tail distributions}
\label{Tail}

In this section, we will state and prove 
the formula for the tail distribution of the sum of independent,
real valued, random variables that satisfy the L\`evy Property.

If one restricts the formula to the case of sums of
independent, identically distributed random variables, one obtains a
formula very similar to
the main result of Hahn and Klass (1997).  
The main differences are that their inequality involves one sided
inequalities, and also that their inequality is more
precise.

This formula also has a strong resemblance to the result of Lata{\l}a.
As we shall show in Section~\ref{Lp}, computing the $L_p$ norm of $U$ is
effectively equivalent to computing $U^*(e^{-p})$.  Then if one
notices that $(1+x)^p$ is very close to $e^{xp}$ for small positive $x$, one
can see that this result and the result of Lata{\l}a are very closely related.
Presumably one could derive Lata{\l}a's result by combining Theorem~\ref{tail} 
with
Theorem~\ref{lp}.  However the technical difficulties are quite tricky, and 
since
Lata{\l}a's proof is elegant, we will not carry out this program here.

\begin{thm}
\label{tail}
Let $(X_n)$ be a sequence of real valued
independent 
random variables satisfying the strong L\`evy property. Define the
functions $F_1(t)$ and $F_2(t)$ to be $0$ if $t > 1$, and if $0 \le
t \le 1$, 
\begin{eqnarray*}
  F_1(t) &=&
    \inf\left\{ \lambda>0 :\prod_n
    \E(t^{X_n^{(\le \ell(t))}/\lambda}) \le t^{-1} \text{ and }
    \prod_n
    \E(t^{-X_n^{(\le \ell(t))}/\lambda}) \le t^{-1} \right\} , \\
  F_2(t) &=&
    \inf\left\{ \lambda>0 :\prod_n \E(t^{X_n^{(\le
    M^*(t))}/\lambda}) \le t^{-1} \text{ and }
    \prod_n \E(t^{-X_n^{(\le
    M^*(t))}/\lambda}) \le t^{-1} \right\} . 
\end{eqnarray*} 
Then 
$$ S^*(t)
\approxsb t F_1(t) \approxsb t F_2(t) ,$$ where the constants
of approximation depend only upon the strong L\`evy constants of
$(X_n)$.
\end{thm}

Let us start with gaining some understanding of
{\em Orlicz spaces}. There is a huge literature on Orlicz spaces,
see for example Lindenstrauss and Tzafriri (1977). 
Suppose that $\Phi:[0,\infty) \to
[0,\infty]$ is a given function (usually convex with $\Phi(0) =
0$). Then the {\em Orlicz norm\/} of a random variable $X$ is defined 
according to 
the formula 
$$ \snormo X_\Phi = \inf\{ \lambda>0 : \E\Phi(\smodo X/\lambda) \le 1
\}.$$ 
We will be concerned with the special functions
$$\Phi_t(x) = 
{\frac{t^{-x} - 1 }{ t^{-1} - 1}} .$$ 
The following is a special case of 
results that appear in Montgomery-Smith (1992).

\begin{lemma}
\label{orlicz}
For any random variable $X$, and for $t \le 1/4$,
we have that 
$$ \snormo X_{\Phi_t} \approx \sup_{0\le x \le 1} {\frac{\log(t) }{
\log(xt)}}
X^*(x) ,$$ with constants of approximation bounded by $2$.
\end{lemma}

\Proof  Suppose first that $\snormo X_{\Phi_t} \le 1$. Then 
$\E\Phi_t(X) \le 1$, which implies that
$$ x t^{-X^*(x)} 
   \le
   \int_0^1 t^{-X^*(y)} \, dy
   \le
   \E(t^{-\smodo X})
   \le t^{-1} ,$$
that is, 
$X^*(x) \le \log(xt)/\log(t)$.

Conversely, suppose that $X^*(x) \le \log(xt)/\log(t)$ for $0\le x \le
1$. Then $$ E\Phi_t(X/2)\le \int_0^1 \Phi_t\left(\frac{\log(xt)}{
2\log(t)}\right) \, dx = {\frac{2t^{-1/2}-1
}{ t^{-1}-1}} \le 1 .$$ \qed

\Proofof{Theorem~\ref{tail}} Let us start with the proof that $S^*(t
 ) \approxsb t F_1(t)$. Since the random  variables $X_n^{(\le \ell(t))}$ 
are independent, we have that $$ F_1(t) =\inf\left\{ \lambda>0 : 
\E(t^{S^{(\le \ell(t))}/\lambda}) \le t^{-1} \text{ and } 
\E(t^{-S^{(\le \ell(t))}/\lambda}) \le t^{-1} \right\} .$$ 
Now we notice that for any random variable $Y$, and $0\le t \le 1$, we
have that$$ 
{\textstyle {\frac{1}{ 2}}} \E(t^{-\smodo Y}) \le
\max\{\E(t^Y),\E(t^{-Y})\} \le 
\E(t^{-\smodo Y}) .$$ Hence $$ F_1(t) \le\inf\left\{ \lambda>0 : 
\E(t^{-\smodo{S^{(\le \ell(t))}/\lambda}}) \le t^{-1} \right\} = 
\snormo{S^{(\le \ell(t))}}_{\Phi_t} ,$$ and 
$$ F_1(t)\ge\inf\left\{ \lambda>0 :\E(t^{-\smodo{S^{(\le
\ell(t))}/\lambda}}) 
\le 2t^{-1} \right\} = \snormo{S^{(\le \ell(t))}}_{\Psi_t} ,$$ where 
$\Psi_t(x) = {\displaystyle{\frac{t^{-x} - 1}{ 2t^{-1} - 1}}}$. However,
we quickly see that for $x \ge 0$ that if $t\le 1/2$ 
then $\Psi_t(x) \ge {\frac{1}{ 3}} \Phi_t(x) \ge \Phi_t(x/3)$, since $\Phi_t$ is a
convex function. Hence $$ F_1(t)
\approx \snormo{S^{(\le \ell(t))}}_{\Phi_t} $$ 
with constants of 
approximation bounded by $3$.

Next, we apply Lemma~\ref{orlicz}, and we see that
$$ F_1(t) 
   \approx 
   \sup_{0\le x \le 1}{\frac{\log(t) }{ \log(xt)}} (S^{(\le \ell(t))})^*(x) .$$ 
Taking $x = t$, we see that the right hand side is bounded below
by ${\frac{1}{ 2}} (S^{(\le\ell(t))})^*(t)$. 
Also, if $t \le x \le 1$,
then $$ {\frac{\log(t) }{ \log(xt)}} (S^{(\le \ell(t))})^*(x)
\le (S^{(\le \ell(t))})^*(t) .$$ 
Further, by Corollary~\ref{klass-nowicki-decr}, there exist
constants $c_1$ and $c_2$, depending only on the L\`evy constants of
$(X_n)$, such that if $0 \le x \le t \le c_1^{-1}$, then 
\begin{eqnarray*}
  {\frac{\log(t) }{ \log(xt)}} (S^{(\le \ell(t))})^*(x)
  &\le& c_2 {\frac{\log(x) }{ \log(xt)}} \bigl((S^{(\le \ell(t))})^*(t)+
        (M^{(\ell(t))})^*(x/2)\bigr) \\
  &\le& c_2 \bigl((S^{(\le \ell(t))})^*(t) + \ell(t)\bigr) .
\end{eqnarray*} 
Now, applying Proposition~\ref{ell-max}, Proposition~\ref{max}, 
and Lemma~\ref{levy} 
part~\ref{levy-ellequiv}, we finally 
obtain the desired result.

To show that $S^*(t) \approxsb t F_2(t)$ is an almost identical
proof. \qed

\section{$L_p$ norms}
\label{Lp}

The main result of this section establishes the
relationship between the $L_p$ norm of sums of random variables and
their tail distributions.

\begin{thm}
\label{lp}
Given $p_0>0$, if
$p \ge p_0$, and $(X_n)$ is a sequence of independent random variables,
then
$$ \snormo U_p \approx U^*(e^{-p}/4) + \snormo
\ell_p   \approx (U^{(\le \ell(e^{-p}/8))})^*(e^{-p}/4) + \snormo
\ell_p ,$$ 
where the constants of approximation depend only upon
$p_0$.
\end{thm}

We should note that we are not able to get universal control over
the constants as $p_0 \to0$, as is shown by
simple examples once one understands that $\snormo Y_p$ converges to
the geometric mean of $\smodo Y$ as $p \to 0$.

Combining this with Corollary~\ref{klass-nowicki-decr},
we immediately obtain the  following result that compares $\snormo S_q$
to $\snormo S_p$. This result extends results of Talagrand, (see
Ledoux and Talagrand (1991, Theorem 6.20), Kwapie\'n and Woyczy\'nski
(1992, Proposition 1.4.2 and comments following it;  see also 
Hitczenko (1994, Proposition 4.1))   and Johnson, Schechtman and Zinn 
(1983).  If this result is specialized to symmetric or positive real
valued random variables, then by considering the cases $p=2$ or $p=1$, 
it implies the inequality of Rosenthal
(1970), including the 
result of Johnson, Schechtman and Zinn (1983) that gives
correct order of the constants as $p \to \infty$.

\begin{thm}
Let $(X_n)$ be a sequence of independent random
variables and let $p_0 > 0$. Then there exist positive constants
$c_1$, $c_2$ 
and $c_3$, depending only upon $p_0$, such that for $q \ge p \ge p_0$
we 
have 
\begin{eqnarray*}
  \snormo U_q
  &\le& c_1 {\frac{q }{ \max\{p,\log(e+q)\}}}
        \bigl(\snormo U_p + M^*(c_2^{-1} e^{-q}) \bigr) + 
        c_1 \normo M_q \\
  &\le& c_3 {\frac{q }{ \max\{p,\log(e+q)\}}}
        \bigl(\snormo U_p + \snormo M_q \bigr). 
\end{eqnarray*}
\end{thm}

Let us proceed with the proofs.  First we need a lemma that allows 
us to deal with the ``large'' parts of $U$, so that they might be 
effectively considered as a sum of disjoint random variables.

\begin{lemma}
\label{disjoint}
Let $(X_n)$ be a sequence of independent random 
variables, and let $0<r<1$.  Then we may express $ U^{(>\ell(r))} = 
\sum_{k=1}^\infty V_k $, where the random variables $V_k$ are
disjoint, and $ V_k^*(t) \le k \ell\bigl(t(k-1)!/r^{k-1}\bigr) $.
\end{lemma}

\Proof  In proving this result, we may suppose without loss of
generality that $X_n = X_n^{(>\ell(r))}$, that is, we may suppose that
$\sum_n\Pr(X_n \ne 0) \le r$.

If $A$ is a finite subset of $\N$, define the event 
$$ E_A = \{ X_n \ne 0\text{ if and only if } n \in A \} .$$ For each
positive integer $k$, let $E_k = \bigcup_{\substack{A \subseteq \N \\ \smodo
A = k}} 
E_A$. Set $V_k = U I_{E_k}$. Notice that if $\smodo A = k$, then 
\begin{eqnarray*}
  \Pr(U I_{E_A} > x) 
  &\le& \sum_{n \in A}
        \Pr(\smodo{X_n} > x/k\text{ and } E_A) \\
  &=& \sum_{n\in A} \Pr(\smodo{X_n} > x/k) \prod_{m \in A \setminus\{n\}}
      \Pr(X_m \ne 0) . 
\end{eqnarray*}
Hence, 
\begin{eqnarray*}
  \Pr(V_k > x) 
  &=& \sum_{\substack{A \subseteq \N \\ \smodo A = k}}\Pr(U I_{E_A} > x) \\
  &\le& \sum_{i_1 < \dots < i_k} \sum_{j=1}^k \Pr(\smodo{X_{i_j}} > x/k) 
        \prod_{\substack{l=1 \\ l\ne j}}^k \Pr(X_{i_l} \ne 0) \\ 
  &\le& {\frac{1}{ k!}} 
        \sum_{i_1} \dots \sum_{i_k}\sum_{j=1}^k \Pr(\smodo{X_{i_j}} > x/k) 
        \prod_{\substack{l=1 \\ l\ne j}}^k \Pr(X_{i_l} \ne 0) \\
  &=& {\frac{k }{ k!}}
      \sum_{i_1} \dots \sum_{i_k}\Pr(\smodo{X_{i_1}} > x/k)
      \prod_{l=2}^k \Pr(X_{i_l} \ne 0) \\
  &=& {\frac{k }{ k!}} 
      \left(\sum_n 
      \Pr(\smodo{X_n} > x/k) \right)\left(\sum_n \Pr(X_n \ne 0)\right)^{k-1} \\
  &\le& {\frac{r^{k-1} }{ (k-1)!}} \Pr(\ell > x/k) . 
\end{eqnarray*} 
\qed

\begin{cor}
\label{disjoint-lp}
Let $(X_n)$ be a sequence of independent random 
variables, let $0<r<1$, and let $0<p<\infty$.  Then
$$ \snormo{U^{(>\ell(r))}}_p \le 2 e^{2^p r/p} \snormo \ell_p .$$
\end{cor}

\Proof  Apply Lemma~\ref{disjoint} to obtain the $V_k$.  
Using the fact that $k \le 2^k$, we obtain that
$$ \snormo{V_k}_p^p
   \le k^p {\frac{r^{k-1}}{ (k-1)!}} \snormo{\ell}_p^p 
   \le 2^p {\frac{(2^p r)^{k-1}}{ (k-1)!}} \snormo{\ell}_p^p .$$
Thus 
$$ \snormo {U^{(>\ell(r))}}_p^p 
   = \sum_{k=1}^\infty \snormo{V_k}_p^p
  \le 2^p \snormo\ell_p^p \sum_{k=1}^\infty {\frac{(2^p r)^{k-1}}{ (k-1)!}}
  = 2^p e^{2^p r} \snormo\ell_p^p .$$ 
\qed

\Proofof{Theorem~\ref{lp}} Applying Proposition~\ref{ell-max}, we see that 
$$ \snormo U_p \ge {\frac{1}{ 2}} \snormo M_p \ge 2^{-1-1/p} \snormo\ell_p
.$$ Also, we have that $$ \snormo U_p^p = \int_0^1 (U^*(t))^p \, dt
\ge 8^{-1} e^{-p} (U^*(e^{-p}/8))^p ,$$ that is, $\normo U_p \ge
8^{-1/p} 
e^{-1} U^*(e^{-p}/8) \ge 8^{-1/p} e^{-1} U^*(e^{-p}/4)$. 
Hence we have shown that there exists a constant $c_1$, depending only
upon $p_0$, such that $$ \snormo U_p \ge c_1^{-1}(U^*(e^{-p}/4) +
\snormo 
\ell_p) .$$ Furthermore, by Proposition~\ref{ell-max}, $$ \Pr(U \ne U^{(\le 
\ell(e^{-p}/8))}) \le \Pr(M > \ell(e^{-p}/8))   \le \Pr(M >
M^*(e^{-p}/8)) \le e^{-p}/8 .$$ Hence $U^*(e^{-p}/8) \ge (U^{(\le
\ell(e^{-p}/8))})^*(e^{-p}/4)$, and so we have shown that there is a
constant $c_2>0$, depending only upon $p_0$, such that $$ \snormo U_p
\ge c_2^{-1}(U^{(\le \ell(e^{-p}/8))})^*                (e^{-p}/4) +
\snormo \ell_p) . $$ Now let us derive the converse inequalities. 
Corollary~\ref{klass-nowicki-decr} tells us that for $t \le e^{-p}/2$ that 
$$ (U^{(\le \ell(e^{-p}/8))})^*(t) \le c_2 {\frac{\log(1/t) }{ p +
\log(2)}}
\bigl((U^{(\le \ell(e^{-p}/8))})^*(e^{-p}/2) + \ell(e^{-p}/8)
\bigr). $$ 
Thus
\begin{eqnarray*}
  \snormo{U^{(\le \ell(e^{-p}/8))}}_p^p
  &\le& \int_0^{e^{-p}/2} ((U^{(\le \ell(e^{-p}/8))})^*(t))^p \, dt 
+(1-e^{-p}/2)(U^{(\le \ell(e^{-p}/8))})^*(e^{-p}/2)  \\ 
  &\le& {\frac{1}{ (p + \log(2))^p}} \int_0^1 (\log(1/t))^p \, dt
        \,\bigl((U^{(\le \ell(e^{-p}/8))})^*(e^{-p}/2) + 
        \ell(e^{-p}/8) \bigr)^p \\
& & \qquad + (U^{(\le \ell(e^{-p}/8))})^*(e^{-p}/2)\\
  &\le& c_3^p\bigl((U^{(\le \ell(e^{-p}/8))})^*(e^{-p}/2) +
       \ell(e^{-p}/8) \bigr)^p ,
\end{eqnarray*} 
where $c_3>0$ depends only upon $p_0$. Furthermore, 
$$ \ell(e^{-p}/8)^p \le 8 e^p \int_0^1 \ell(t)^p \, dt = 8 e^p
\snormo\ell_p^p.$$ Hence, applying Corollary~\ref{disjoint-lp}, and the
(quasi-)triangle 
inequality for $L_p$, we deduce that
 there exists a constant $c_4$, depending only upon $p_0$, such that 
$$ \snormo U_p \le c_4(U^{(\le \ell(e^{-p}/8))})^*  (e^{-p}/2)
 + \snormo \ell_p) . $$ Finally the result follows by noticing that 
$$ (U^{(\le \ell(e^{-p}/8))})^*(e^{-p}/2) \le   (U^{(\le
\ell(e^{-p}/8))})^*
(e^{-p}/4) ,$$ and also, by an argument similar to one presented
above, 
that $$ (U^{(\le \ell(e^{-p}/8))})^*(e^{-p}/2) \le U^*(3e^{-p}/8) \le
U^*(e^{-p}/4) .$$ \qed

\section{Rearrangement invariant spaces}
\label{Ri}

Rearrangement invariant spaces are studied in much 
of the
literature, see for example Lindenstrauss and Tzafriri (1977).
However, we will work with a definition that is a little
less restrictive.
A {\em rearrangement invariant space\/} on the
random variables is a quasi-normed Banach space $\cL$ of
random variables such that $1 \in \cL$, and if $X^* \le Y^*$ and $Y
\in \cL$, then $X \in \cL$ and $\normo X_\cL \le \normo Y_\cL$.
Obviously the spaces $L_p$ for $0<p\le \infty$ are rearrangement
invariant spaces.

Given a rearrangement invariant space $\cL$, we
define the {\em quasi-constant\/} of $\cL$ to be the least constant $K>0$ such
that $\normo{X+Y}_\cL \le K(\normo X_\cL + \normo Y_\cL)$ for
 all $X,Y \in \cL$.  Notice that
 if $X^*(2t) \le Y^*(t)$, and $Y \in \cL$, 
then $X$ may be written as the sum of two disjoint 
random variables $Y_1$ and $Y_2$ with $Y_1^*(t),Y_2^*(t) \le Y^*(t)$,
and hence $\normo X_\cL \le 2K \normo Y_\cL$.

Given two rearrangement invariant spaces $\cL$ and $\cM$, we will say
that $\cL$ embeds into $\cM$ if there is a positive constant $c$ such that if
$X \in \cL$, then $X \in \cM$ and $\normo X_\cM \le c \normo X_\cL$.
We will call the least such $c$ the {\em embedding constant\/} of
$\cL$ 
into $\cM$.

\begin{thm}
Let $p_0>0$, and let $\cL$ be a rearrangement invariant
space such that $\cL$ embeds into $L_p$, and $L_q$ embeds into $\cL$, 
where $q \ge p \ge p_0$. Then there is a positive constant $c$, depending 
only upon the quasi-constant of $\cL$, the embedding constants, $p_0$
and 
$q/p$, such that 
for any sequence of independent random variables $(X_n)$
$$ c^{-1}(\normo U_p + \normo\ell_\cL) \le \normo
U_\cL  \le 
c(\normo U_p + \normo\ell_\cL) .$$
\end{thm}

\Proof   Let us first obtain the left hand side
inequality. It follows by hypothesis that $\normo U_\cL \ge c_1^{-1}
\normo U_p$, where $c_1$ is the embedding constant of $\cL$ into
$L_p$. Furthermore, $U \ge {\frac{1}{ 2}} M$, and by Proposition~2.1, 
$\ell(t) \le M^*(2t)$. Hence $ \normo U_\cL \ge (4K)^{-1}
\normo\ell_\cL$, where $K$ is the quasi-constant of $\cL$. 

Now let us obtain the right hand inequality. 
By Corollary~\ref{klass-nowicki-decr}, we
have that 
there is a universal positive $c_2$ for $0 \le t \le 1$ 
$$ U^*(t) I_{0\le t \le 2^{-2q/p}} \le c_2 {\frac{2q }{ p}}
(U^*(t^{p/2q}) + 
M^*(t/2)). $$ 
Now $U^*(t) I_{0\le t \le 2^{-2q/p}} \ge U^*(2^{2q/p}
t)$, and hence $$ \normo U_{\cL} \le (2K)^{\lceil 2q/p \rceil} c_2 K
{\frac{2q}{ p}} (\snormo{t\mapsto U^*(t^{p/2q})}_\cL + K \normo M_\cL) .$$
To finish the 
proof, suppose that $\normo U_p = \lambda$.  Then it is easily seen that
$U^*(t) \le \lambda t^{-1/p}$.  Thus, if $c_3$ is the embedding
constant of 
$L_q$ into $\cL$, then 
\begin{eqnarray*}
  \snormo{t\mapsto U^*(t^{p/2q})}_\cL   
  &\le& c_3 \snormo{t\mapsto U^*(t^{p/2q})}_q \\
  &=& c_3 \left(\int_0^1 (U^*(t^{p/2q}))^q \, dt \right)^{1/q} \\  
  &\le& c_3 \lambda \left(\int_0^1 t^{-1/2} \, dt\right)^{1/q} \\  
  &=&   2^{1/q} c_3 \lambda . 
\end{eqnarray*} 
\qed

\end{document}